\documentclass[11pt]{article}

\usepackage{amsfonts}

\newenvironment{debut}[1]%
{\begin{quote} \bf #1 \it}%
{\end{quote}}

{\begin{debut} {R\'esum\'e~: }}%
{\end{debut}}

\newenvironment{abs}%
{\begin{debut} {Abstract:}}%
{\end{debut}}

\title{On the number of return words in infinite words with complexity
$2n+1$}
\author{Laurent Vuillon\footnote{LIAFA, Universit\'e Paris 7, case 7014,
2 place Jussieu, F-75251 Paris Cedex 05, France,
 vuillon@liafa.jussieu.fr}}
\date{ }

\def\cqfd{\nobreak\nopagebreak\rule{0pt}{0pt}\nobreak\hfill\nobreak
          \rule{.5em}{.5em}}
\newtheorem{remq}{Remarque}

\newtheorem{theo}{Theorem}
\newtheorem{prop}{Proposition}

\newtheorem{lem}{Lemma}
\newtheorem{cor}{Corollary}

\newcommand{\GN}{\mathbb{N}}

\newcommand{\GZ}{\mathbb{Z}}
\newcommand{\GR}{\mathbb{R}}

\def\ind{_{n \in \GN}}

\begin{document}
\maketitle

\begin{abs}
In this article, we count the number of return words in some
infinite words with complexity $2n+1$. We  also consider  some infinite
words
given by  codings of  rotation and interval exchange transformations on $k$
intervals. We prove that the number of return words over a given word $w$
for these infinite words
is exactly $k.$
\medskip


\end{abs}

\section{Introduction}

The starting point of this article is the following
question: {\em
Does there exist a way to characterize  by return words the infinite word
with
complexity $2n+1$~?
}
We first give the definition of return words for a recurrent infinite word
in
a finite
alphabet (a recurrent word is  infinite word such that each word appears
infinitely many
times). Considering
each occurrence of a word $w$ in a recurrent infinite word $U,$
we define the set of return words over $w$ to be the set of all distinct
words
beginning with an occurrence of $w$ and ending exactly before the next
occurrence of $w$ in the infinite sequence. This mathematical tool was
introduced independently  by Durand, Holton and Zamboni in order to
study
primitive substitutive sequences (see \cite{D,DT,Ho}).
This notion is quite natural and can be seen as a symbolic version
of the first return map for a dynamical system. Recently,
many developments of the notion    of return words
have been given. For example, Allouche, Davinson, Queffelec and Zamboni
study the transcendence of Sturmian or morphic continued
 fractions and the main tool is to show, using return words, that
arbitrary long prefix are
 ``almost squares'' (see \cite{ADQZ}).
We can also use return words to study low complexity infinite words.
For example,  the  author shows that an infinite word is Sturmian
if and only if for each  word appearing in the infinite word, the
cardinality of the set of return words over $w$ is exactly two (see
\cite{V2}).
Recall that  Sturmian infinite words are aperiodic  words with complexity
$p(n)=n+1$
for all $n$ (the complexity function $p(n)$ counts the number of
distinct
factors of length $n$ in the infinite word) (see \cite{B,BS,H}).
 Fagnot and Vuillon show a
generalization of
the notion of balanced property for Sturmian words and the
proof is based on return words and combinatorics on words (see
\cite{FV}).
Cassaigne uses this tool to investigate a Rauzy's conjecture (see
\cite{C}).

 Our propose
is to compute the set of return words for a class of infinite words with
complexity
$p(n)=2n+1$ for all $n$ (see \cite{AR,BV,CFZ,CMR,DJP,F,JV}).

The following question is thus natural: is it true that for these
infinite words,
 the cardinality of the  set of return words over
a factor $w$ is always three ?

In this direction, Justin and Vuillon \cite{JV} show that
  Arnoux-Rauzy infinite words \cite{AR}, which are
infinite words with complexity $2n+1$ (see \cite{AR,CFZ,CMR,DJP}),
 have property  $R_3.$
In the sequel, we say that an infinite word have the property $R_n$
if the cardinality of the  set of return words over
each factor $w$ is exactly $n$.
They give also the structure of return words in the context of
Episturmian words \cite{DJP,JV}.

In the opposite, the work of Ferenczi \cite{F} presents a nice substitutive
infinite word with
complexity  $p(n)=2n+1$ given by the Chacon transformation.
It is not difficult to show that the cardinality of the set of return
words may be
upper than three for this infinite word. This is the first example
of infinite word with complexity 2n+1 without the property  $R_3.$

Nevertheless, we find two other classes of infinite words with complexity
$p(n)=2n+1$ having the property  $R_3.$
The first one is the coding of a rotation $\alpha$ in the unit circle with
three intervals (with rationaly independent lengthes) and one of these
intervals
of length $\alpha.$ The second one
is a generalization of this class namely the interval exchange
transformation
on three intervals.

The structure of the article is the following.
First,
 we recall some basic definitions in combinatorics of words.
Secondly, we present the infinite words given by codings of rotation.
Then, we show that an infinite word given by a regular  interval
exchange with  three intervals (resp. $n$ intervals) have the property
$R_3$ (resp. $R_n$). We compute also the length of each return word over a
given finite
word $w.$ At last, we show that  the coding of a rotation in the
unit circle with
three  (resp. $n$ intervals) intervals with rationaly independent
lengthes as property $R_3$ (resp. $R_n$)

More generally, we guess that the class of infinite word with complexity
 $p(n)=2n+1$
and
exactly three return words over each factor is  the class of infinite words
given by a ``self-induced  discrete dynamical system''.
In particular, we may extend these results to bounded remainder sets
introduced by Rauzy and Ferenczi.

\section{Basic definitions and examples}

Let $\mathcal{A}=\{0,1\}$ be a binary alphabet. We denote by
$\mathcal{A}^*$ the set of finite
words on  $\mathcal{A}$ and
by  $\mathcal{A}^{\GN}$ the one-sided infinite word $U.$
A       word $w$ is a {\em factor} of a word $x \in \mathcal{A}^*$ if
there
exist  some words $u,v \in \mathcal{A}^* $ such that $x=uwv.$
An       infinite word $U$ is called {\em recurrent} if
every factor of $U$ appears infinitely many times in $U.$
For a finite word $w=w_1w_2 \cdots w_n,$ the length of $w$  is denoted
 by $|w|$ and is equal to $n.$ The set of factors of $U$ with length
 $n$ is denoted by  $L_n(U).$ The language $L(U)=\cup_n L_n(U)$ is the
set of factors of $U.$
For two finite  words $w$ and $u,$ the number of occurrences of $w$ on
$u$ is denoted
by $|u|_w$
and $|u|_w=\{i|\  0 \leq i \leq |u|-|w| \mbox{ s.t. } u_iu_{i+1} \cdots
u_{i+|w|-1} = w_1w_2 \cdots
w_{|w|}\}.$

The  {\em position set} ($i(U,w)=\{i_1,i_2,\cdots, i_k, \cdots\}  $) of
the word $w$ is a set of integers $i(U,w)=\{i_1,i_2,\cdots, i_k,
\cdots\}$ where  $i_k$ represents the
position of
the first letter of  the
$k$-th occurrence
of the word $w$ in the infinite word $U.$
In a more formal way,
 $i_k \in i(U,w)$ if and only if  $U_{i_k}U_{i_k+1}\cdots
U_{i_k+|w|-1}=w$ and $|U_1 \cdots U_{i_k+|w|-1}|_w=k.$
Since the infinite word $U$ is recurrent
the set $i(U,w)$ is infinite.
For a recurrent          word $U,$  the set of  {\em return words over
$w$}
is the set (denoted by $ \mathcal{H}_{U,w}$)
of all distinct words with the following form:
$$U_{i_k}U_{i_{k}+1}\cdots U_{i_{k+1}-1}$$
for all $ k \in \GN, k>0.$
This definition is best understood on an example.
Let $U_1=(0100100001)^\omega$ be  an  infinite word on the alphabet
$\mathcal{A}.$
By definition, the set of return words over $01$  is
$ \mathcal{H}_{U_1,01}=\{010,01000,01\}.$
Indeed, the infinite word $U_1$ can be written
$$(\underline{0}10\underline{0}1000\underline{0}1\underline{0}
10\underline{0}1000\underline{0}1)^\omega$$
where $\underline{~}$ denotes  the position of the first letter for each
occurrence of the word 01.
We note that between two consecutive occurrences of 01, we have  the
three  possible words  $0,000,
\epsilon.$
Thus 010,01000 and 01 are elements of $ \mathcal{H}_{U_1,01}.$
We say that an infinite word have the property $R_n$
if the cardinality of the  set of return words over
each factor $w$ is exactly $n$ (i.e. $Card(\mathcal{H}_{U,w})=2$ $\ \forall
w \in L(U))$.

The {\em complexity function}  $p: \GN  \rightarrow \GN$ of an infinite
word
$U$
 counts the number of  distinct factors of $U$ of given length:
$$p(n)=Card\{w| \ w \in L_n(U) \}.$$

\section{A negative result}

Let consider the following substitution
extensively studied by Ferenczi (see \cite{F})  which is a recoding
of the Chacon substitution:

$$\sigma(1)=12,$$
$$\sigma(2)=312,$$
$$\sigma(3)=3312.$$

The fix point $\sigma(x)=x$ of the substitution begins by
$$x=1231233121231233123312123121231233121231231233 \cdots$$

It is easy to check with few terms of the fix point $x$ that
 the number of return words over $23$ is upper
than 3, indeed:

$$\mathcal{H}_{23}=\{231,2331,23121,233121\}.$$

\section{Codings of rotation}
The aim of this section is to introduce  codings of irrational
rotation
on the unit circle.
For $p \geq 2$, let $F= \{\beta_0<\beta_1<\ldots<\beta_{p-1}\}$ be a set
 of $p$
consecutive points of the unit circle (identified in all that
 follows with $[0,1[$ or with the unidimensional torus $\GR/\GZ$) and
let $\beta_p=
\beta_0$.
Let $\alpha$ be an irrational number in $]0,1[$ and let us consider
the {\em positive orbit} of a point $x$ of the unit circle under the
rotation by angle $\alpha$, i.e.,     the set
of points $\{ \{\alpha n +x\} ,\ n
\in \GN\}$.

The {\em coding}  of the  orbit of $x$ under the rotation by angle
$\alpha
$
 with respect to  the partition
$\{[\beta_0,\beta_1[,[\beta_1,\beta_2[,\ldots,[\beta_{p-1},\beta_{p}[\}$
is the infinite word $u$
defined on the
finite alphabet
$\Sigma= \{0,\ldots,{p-1}\}$
as
follows:
$$u_n= k \Leftrightarrow \{x+n\alpha\} \in [\beta_k,\beta_{k+1}[,\
\mbox{ for }
\ 0\leq k\leq p-1.$$
 {\em  A coding of the rotation $R$}
means  the  coding of the  orbit of  a point $x$ of the unit circle
under the rotation $R$
 with respect to  a finite partition of the unit circle consisting of
left-closed and right-open intervals.

\subsection{Factors}

With the above notation,   consider a coding $u$ of the
orbit of a point $x$ under the rotation by angle $\alpha$
 with respect to the partition
$\{[\beta_0,\beta_1[,[\beta_1,\beta_2[,\ldots,[\beta_{p-1},\beta_{p}[\}.$
Let $I_k= [\beta_k,\beta_{k+1}[$ and let $R$ denotes the rotation by
angle
 $\alpha$.
\begin{lem}
 $w_1\cdots w_n$
defined on the alphabet $\Sigma=\{0,1,\ldots,p-1\}$ is a
 factor of the infinite word $u$
if and only if
$I(w_1,\ldots,w_n)= \bigcap_{j= 0}^{n-1}R^{-j}(I_{w_{j+1}})  \neq
\emptyset.$
\end{lem}

\paragraph{Proof}
A  finite word $w_1\cdots w_n$
defined on the alphabet $\Sigma=\{0,1,\ldots,p-1\}$ is a
 factor of the infinite word $u$
if and only if there exists an integer $k$ such that
 $$\{x+k\alpha\} \in
I(w_1,\ldots,w_n)= \bigcap_{j= 0}^{n-1}R^{-j}(I_{w_{j+1}}).$$
As $\alpha$ is irrational, the sequence $(\{x+n\alpha\})\ind$ is
dense in the unit circle, which implies that $w_1w_2\ldots w_n$ is a
factor
of  $u$
if and only if $I(w_1,\ldots,w_n) \neq \emptyset$.
In particular, the set of factors of a coding does not depend on the
initial point $x$ of this coding.
Furthermore, the connected components of these
sets are bounded by the points
$$\{k(1-\alpha)+ \beta_i\} \mbox{, for } 0 \leq k \leq n-1,\ 0\leq i
\leq p-1.$$

These sets consist of  finite unions of intervals.
More precisely, if for every $k$, $\beta_{k+1}-\beta_k \leq
\sup(\alpha,1-\alpha)$, then these
 sets
 are connected;  if there exists $K$ such that
$\beta_{K+1}-\beta_K > \sup(\alpha,1-\alpha)$, then the sets are
connected
except for $w_1\ldots w_n$
of the form $a_{K}^n$ (see \cite{albe}) (the notation $a_K^n$
denotes the word of length $n$  obtained by successive concatenations of
the letter $a_K$).  Let us note
that there exists at most one integer  $K$ satisfying
$\beta_{K+1}-\beta_K > \sup(\alpha,1-\alpha)$. \cqfd

\section{Interval exchange transformations}

This section deals with interval exchange
transformations. This object is a natural generalization of the codings
of
rotation $\alpha$ where $\alpha$ is the length of one of the intervals.
An interval exchange transformation
is a piecewise affine transformation which maps a partition  into
intervals
of the space
to another
partition according to a permutation. This transformation could be
more complicated than a rotation. Indeed, irrational rotations
are uniquely ergodic, unless a class of interval exchange
transformations  is non uniquely ergodic. In the sequel, we use the
notations of Keane and Rauzy (see \cite{K,Ri})

 An  interval exchange of $k$ intervals is defined by a vector
$(\lambda_1,\lambda_2, \cdots, \lambda_k)$ in $\GR^k$
with strictly positive coordinates and $\sum_{i=1}^k \lambda_i = 1$,
 and a permutation $\sigma$ of
the set $\{1,2,\cdots,k\}$. We set the partition of
$[0,1[$
in $k$ intervals $$X_i= \left [\sum_{j<i} \lambda_j ,\sum_{j \leq i}
  \lambda_j \right [ \mbox{ for } 1 \leq i \leq k$$
($X_i$ has the length $\lambda_i$).

The interval exchange transformation associated to the ordered pair
$(\sigma, \lambda)$ is the transformation from $[0,1[$ to itself,
defined as a piecewise affine transformation which maps the partition
$(X_1,X_2, \cdots, X_k)$ to the
partition $(X_{\sigma(1)},X_{\sigma(2)}, \cdots, X_{\sigma(k)}).$

The transformation $T$ maps the point $x \in X_i$
to the point $$T(x)= x+a_i$$
where $$a_i=\sum_{k<\sigma^{-1}(i)}   \lambda_{\sigma(k)}-
\sum_{k<i} \lambda_{k}.$$

Now, we construct an  infinite word  with values in the alphabet
$\{1,2, \cdots,k\}$
associated with a couple $(\sigma, \lambda)$
by coding the positive  orbit of a point $x$ by the transformation $T$
according to the partition  $(X_1,X_2, \cdots, X_k).$
The  infinite word $U(x)$ is given by
$$U(x)_n= \mathcal{I}(T^n(x)),~\forall n \in \GN$$ where $\mathcal{I}(y)=i$
if $y \in X_i.$

In order to define a transformation such that  each orbit
is  dense in $[0,1[$, we add the following property:
an interval exchange transformation is called regular if
for all  $0=a_1 < a_2 < \cdots < a_{k+1}=1$
of the intervals $X_i=[a_i,a_{i+1}[$ with $i \in \{ 1,\cdots,k \},$
we have $$T^n a_i = a_j,~ i \mbox{ and } j \in \{2,3, \cdots, k\}, n \in
\GZ$$
implies $n=0, i=j.$

M. Keane shows the following result \cite{K}:
\begin{theo}
An interval exchange transformation $T$ is regular if and only if for
all point $x$ in  $[0,1[$ the orbit $\cup_{n \in \GZ} \{T^n x\}$ of
$x$ by $T$ is dense in  $[0,1[.$
\end{theo}

A necessary condition to have a regular interval exchange
transformation
is to take an irreducible permutation.
A permutation is called irreducible if no subset of
$\{1,2,\cdots,k\}$ is invariant by the permutation,
i.e.,  $\sigma(\{1,2, \cdots, \ell\}) \neq \{1,2, \cdots, \ell\}$
for all $\ell<k.$

Now we can state the main theorem of this section:
\begin{theo}
 Let a  regular  interval exchange
transformation $T$ with $k$ intervals. The  infinite word
$U(x)$
associated with $T$ has property $R_k$.
\end{theo}

In order to show this theorem, we first  present a method to construct
all the words of length $n$ factors of the infinite word $U(x)$ and secondly
a way to construct return words by using self-induction.

\subsection{Factors}

The construction of all the factors is not so far from
the one used for rotations on the unit circle. The main tool
is to consider the negative orbit of all the
 endpoints $0=a_1 < a_2 < \cdots < a_{k+1} =1$
of the intervals $X_i=[a_i,a_{i+1}[$ with $i \in \{1, \cdots,k\}.$
The set of endpoints $\{a_i |  \ 1\leq i \leq k+1\}$ is called  $X^{(1)}.$

A word $w$ of length $n$
is a factor of the infinite word $U(x)$ if and only if
there exists an interval $I_w \subset  [0,1[$ and a point $y$ in $I_w$
such that
the word $w$ is the following
concatenation of letters: $$\mathcal{I}(x)\mathcal{I}(T(x))\cdots
\mathcal{I}(T^{n-1})=w.$$

Using this fact,   the number of factors of length one is
exactly the number of  intervals $I_1=X_1, I_2=X_2, \cdots, I_k=X_k$
associated with the letters of the alphabet $\{1,2, \cdots, k\}.$
\\

\begin{prop}
Let a  regular  interval exchange
transformation $T$ with $k$ intervals. The infinite  word
$U(x)$
associated with $T$ has complexity $p(n)=n(k-1)+1, \forall n \in \GN$.
\end{prop}

\paragraph{Proof}
The factors $w=w_1w_2$ of length 2 are given by the set
$\{\mathcal{I}(x)\mathcal{I}(T(x)), x \in [0,1[\}$
or equivalently by the intervals $I_{w_1w_2}= I_{w_1}\cap T^{-1}
I_{w_{2}}.$

As the transformation $T$ is a piecewise affine transformation,
it is sufficient to find all the endpoints of the intervals
$I_w$ with $|w|=2.$
That is the positions of the endpoints of all intervals
with form  $I_{w_1}\cap T^{-1} I_{w_{2}}.$
These endpoints
are given by the ordered set
$$X^{(2)}=\{0=b_1 < b_2 < \cdots < b_{2k}=1\}= X^{(1)} \cup
T^{-1}X^{(1)}.$$
Remark that the number of points in the union is $2k$
because  the points $0$ and $1$ are elements of the intersection
of the two sets  $X^{(1)}$ and $T^{-1}X^{(1)}.$

The next step is to prove that each $I_w$ is connected.
By definition, all the intervals associated with words
of length $1$ are connected. The interval $I_{m_1m_2}= I_{m_1}\cap
T^{-1}
I_{m_{2}}$
is  connected  because it is the intersection of two intervals in
$[0,1[.$

In other words, for each word $w$
with $|w|=2$ factor of the infinite word $U(x),$
there exists an interval $I_w=[b_\ell,b_{\ell+1}[$ where
 $b_\ell$ and $b_{\ell+1}$ are two consecutive points in the ordered set
 $X^{(2)}.$
\\

By induction, let $w=w_1w_2\cdots w_n$ a factor of length $n,$
then there exists a point $x$ such that
$$w=\mathcal{I}(x)\mathcal{I}(T(x))\cdots \mathcal{I}(T^{n-1}(x))$$
or equivalently
$$I_{w_1w_2\cdots w_n}= \cap_{i=0 \cdots n-1} T^{-i} I_{w_{i+1}}.$$
That is each endpoints of the interval associated
with a word of length $n$ is an element of the ordered set
$$X^{(n)}=\{0=b_1 < b_2 < \cdots < b_{n(k-
1)+2}=1\}=
X^{(1)} \cup T^{-1}X^{(1)} \cup \cdots \cup T^{-n+1}X^{(1)}.$$
Remark that the number of points in the union is $n(k-1)+2$
because at each step, we add $k-1$ new points and for $n=1$
we have $k+1$ points in the partition.

Furthermore, the intervals associated with factors
of length $n$ are connected. Indeed, the intervals are the
intersection
of connected intervals in $[0,1[.$

In other words, for each word $w$
with $|w|=n$ factor of the infinite word $U(x),$
there exists an interval $I_w=[b_\ell,b_{\ell+1}[$ where
 $b_\ell$ and $b_{\ell+1}$ are two consecutive points in the ordered set
 $X^{(n)}.$

With this construction, we find $n(k-1)+2$ points
for the endpoints in the partition   $X^{(n)}.$
We have $n(k-1)+1$ intervals in the partition and then the complexity
function
for a regular interval exchange in $k$ intervals is equal to
$p(n,k)=n(k-1)+1, \forall (n,k) \in \GN^2.$ \cqfd

Thus an infinite word associated with
a  regular  interval exchange with  two intervals
is nothing but a Sturmian infinite word (i.e. an infinite word with
complexity
$p(n,2)=n+1$ for all $n).$
And an infinite word associated with
a  regular  interval exchange with  three intervals
is  an infinite word with complexity
$p(n,3)=2n+1$ for all $n.$

\subsection{self-induction}

Now, we focus on the construction of return words
associated with a word $w$ factor of the infinite word.
\medskip

Theorem 2 {\em Let a  regular  interval exchange
transformation $T$ with $k$ intervals. The infinite  word
$U(x)$
associated with $T$ has property $R_k$.}
\medskip

By the previous construction, we find a unique
connected interval $I_w$ associated with the word $w.$
The main tool is to  study the first return map in the
adherence of the  interval
$I_w.$ This method is used by Rauzy to give a generalized continued
fraction
algorithm \cite{Ri,Rii}.

\paragraph{Proof}
Let $r(y)=\inf_{t \geq 0} \{ T^{-t}y \in \overline{I_w} \}$ be the negative
first
return time in
the
interval $ \overline{I_w}.$
Keane shows that the points $\{T^{-r(y)}y | y \in  \overline{I_w}, k(y)<
\infty\}$
give a partition of  the interval $ \overline{I_w}$ in exactly $k$
intervals
$I_{p_1},I_{p_2}, \cdots, I_{p_k}$ (see \cite{K}).
Indeed, the endpoints of these intervals are given by the
first time that the negative orbit of the points in
$X^{(1)}=\{0=a_1 < a_2 < \cdots < a_{k+1}=1\}$ falls in $
\overline{I_w}.$

There are $k+1$ points in $X^{(1)}.$
This shows that for a general interval $[\alpha,\beta]$,
 the number of induced points in $]\alpha,\beta[$ is $k+1$
 and that
the number of induced intervals
is $k+2.$
As the endpoints of $ \overline{I_w}$ are both in the negative orbit
of two different points in $X^{(1)}$.
More precisely, as the interval exchange transformation is regular,
we have  $$T^n a_i = a_j,~ i \mbox{ and } j \in \{2,3, \cdots, k\}, n \in
\GZ$$
implies $n=0, i=j.$
Thus $$\alpha=T^{n_{a_0}}a_0, \beta=T^{n_{a_{k+1}}}a_{k+1}.$$

  That is the number of induced
 points in the interior of $\overline{I_w}$ is $k-1.$
Thus
the number of induced intervals
is $k.$ Such intervals with endpoints on the negative orbit of $a_0$ and
$a_{k+1}$
 are called acceptable intervals by Rauzy (see \cite{Ri}).
Indeed, for interval exchange transformation with $k$ intervals,
  the induced transformation
 on an acceptable intervals
is also an interval exchange transformation with $k$ intervals.

 By construction, for each $p_i$ there exists $y \in I_w$  and $t$ such that
the word
$p_i=\mathcal{I}(y)\mathcal{I}(T(y))\cdots \mathcal{I}(T^{t}(y)) \cdots
\mathcal{I}(T^{t+|w|-1}(y)).$
Furthermore,
 the prefix of length $|w|$ of $p_i$ is  exactly $w$ (indeed $y \in
I_w$) and
the suffix of length $|w|$ of $p_i$ is  exactly $w$ (indeed $t$
is the positive first return time of $y$ in the interval $I_w$
defined by $r^+(y)=\inf_{t \geq 0} \{ T^{t}y \in \overline{I_w} \}$).
In other words, the number of return words over $w$ is exactly the
number
of induced intervals. By construction $I_w$ is an  acceptable interval.
Thus the number of return words over $w$ is exactly $k.$
That is, the infinite word associated with regular interval exchange
transformation
on $k$ intervals has property $R_k.$ \cqfd

\begin{cor}
The length of the i-th return word of $w$ associated with
the i-th induce interval $p_i$ is $|w|+k_i+k'_i$  where $k_i$
is the given by the smallest $k_i \geq 0$ such that
 $T^{-k_i}\ell_i \in I_w,$  where $k_i$ is the given by the smallest
$k'_i>0$ such that
 $T^{k'_i}\ell_i \in I_w$ and where $\ell_i$
is the left endpoint of $p_i.$
\end{cor}

\paragraph{Proof} By construction, the left endpoint of $p_1$ is $\alpha.$
Each orbit of
 the points in
$p_1$ according to the induce transformation is a translation of the orbit
of  $\alpha.$
Then the length of the return word associated with $p_1$ is
$|w|+k'_1$ where $k'_1$ is  the smallest $k'_1>0$ such that
 $T^{k'_1}\alpha \in I_w.$ In other words, the
 length of the return word associated with the interval
$p_1$ is exactly the sum of  the positive first
return time $k'_1$ of the point $\alpha$ in $I_w$ and  the length of $w.$

 By construction,  each interval
 $p_i$ for $i=1,\cdots k$ has $T^{-k_i}\ell_i$ for left endpoint.
Furthermore, each orbit of the points in
$p_i$ according to the induce transformation is a translation of the orbit
of  $T^{-k_i}\ell_i.$
Then the length of the return word associated with $p_i$ is
$|w|+k_i+k'_i$ where $k_i$ is  the smallest $k_i>0$ such that
 $T^{-k_i}\ell_i \in I_w$ and where $k'_i$ is  the smallest $k'_i>0$ such
that
 $T^{k'_i}\ell_i \in I_w.$ In other words, the
 length of the return word associated with the interval
$p_i$ is exactly given by the sum of the first time  that $T^{-k_i}\ell_i$
is mapped on the point $\ell_i$ ($\ell_i$ is element of the partition
$X^{(1)}$) and the positive
 first
return time $k'_i$ of the point $\ell_i$ in $I_w$ and  the length of $w.$
\cqfd

\section{Return words for  codings of rotation}

As the codings of rotation are a special case of interval exchange
transformations,
we have the following corollary:

\begin{cor}
 Let $T$ be a coding of rotations on  $k$ intervals with rationaly
independent lengthes and by a rotation of length $\alpha$ equal to the
length
of one of the intervals. The infinite word $U(x)$
associated with $T$ has property $R_k$.
\end{cor}

\paragraph{Proof}
As the rotation is defined on the circle, we can translate all the
intervals
in order to have the first interval with length equal to $\alpha.$
The coding of rotation of angle $\alpha$ on $k$ intervals with the first
interval
of length
$\alpha$  is equivalent
to an interval exchange transformation with either the permutation
$(2 3 \cdots k 1)$ if $\alpha>0$ or the permutation $(k 1 2 3 \cdots
k-1)$
$\alpha<0.$
As  the lengthes of the intervals are rationaly independent in the coding of
rotation
then the associated
interval exchange transformation is regular.
Consequently, the infinite word $U(x)$
associated with $T$ has property $R_k$.
\cqfd


\begin{thebibliography}{99}


\bibitem{albe}
P. Alessandri, V. Berth\'e {\em Three distance
theorems and combinatorics
on words},
Enseig. Math. {\bf 44} (1998), 103--132.

\bibitem{ADQZ} J.-P. Allouche, J. L. Davison, M. Queff\'elec and
  L. Q. Zamboni {\em Transcendence of Sturmian or morphic continued
  fractions,}
preprint 1999.

\bibitem{AR} P. Arnoux et G. Rauzy,  {\em Repr\'esentation
g\'eom\'etrique
de suites de complexit\'e 2n+1,}
Bull. Soc. Math. France. {\bf 119} (1991), 199--215.



\bibitem{BS} J. Berstel and  P. S\'e\'ebold ,  {\em Sturmian words,}
In M. Lothaire, Algebraic combinatorics on Words.  (2000). To appear.

\bibitem{B} V. Berth\'e, {\em Fr\'equences des facteurs des suites
sturmiennes}, Theoret. Comput. Sci. {\bf 165} (1996), 295--309.

\bibitem{BV} V. Berth\'e and L. Vuillon, {\em A two-dimensional
generalization
 of Sturmian sequences :
tilings and rotations}, to appear in Discr. Math..

\bibitem{C} J. Cassaigne, {\em Ideas for a proof of Rauzy's conjecture
on the
 recurrence functions of infinite words,} talk in Rouen ``words 1999''.

\bibitem{CFZ} J. Cassaigne, S. Ferenczi and L.Q. Zamboni {\em Imbalances in
Arnoux-Rauzy Sequences,} Preprint IML 1999.


\bibitem{CMR} M.G. Castelli, F. Mignosi
 and A. Restivo, {\em Fine and Wilf's theorem for three periods and a
generalization
of Sturmian words,} Theoret. Comp. Sci., to appear.




\bibitem{D} F. Durand,  {\em A characterization of substitutive
sequences using return
words,}
Discrete Math. {\bf 179} (1998), 89--101.

\bibitem{DT} F. Durand,  {\em Contributions  \`a l'\'etude des suites et
syst\`emes dynamiques substitutifs,}
 Ph.D. Thesis, Universit\'e de la M\'editerran\'ee
(Aix-Marseille II) 1996.


\bibitem{DJP} X. Droubay, J. Justin and G. Pirillo, {\em
Episturmian words and some constructions of de Luca and Rauzy},
Theoret. Comp. Sci., to appear.

\bibitem{FV} I. Fagnot and L. Vuillon {\em Generalized balances in
    Sturmian words} Preprint 02/2000 LIAFA.

\bibitem{F} S. Ferenczi {\em Les transformations de Chacon: combinatoire,
structure g\'eom\'trique, lien avec les syst\`emes de complexit\'e 2n+1},
Bull. Soc. Math. France, {\bf 123} n0. 2 (1995), 271--292.

\bibitem{H}  G. A. Hedlund, Morse {\em Symbolic dynamics II:
 Sturmian trajectories}, Amer. J. Math. {\bf 62} (1940), 1--42.

\bibitem{Ho} C. Holton, L. Q.  Zamboni
{\em Descendants of primitive
substitutions,} Theory Comput. Systems 32, (1999) 133--157.

\bibitem{JV} J. Justin and L. Vuillon {\em Return words in Sturmian and
episturmian words,}
Preprint 10/2000 LIAFA.

\bibitem{K} M. Keane {\em Interval exchange transformations}
Math. Zeit. {\bf 141} (1975), 25--31.

\bibitem{Ri} G. Rauzy, {\em Echange d'intervalles et transformations
    induites,} Acta. Arith. {\bf 34} (1979), 315--328.

 \bibitem{Rii} G. Rauzy, {\em Une g\'en\'eralisation du d\'eveloppement
 en fraction continue,} Semin. Delange-Pisot-Poitou, 18e Ann\'ee {\bf
 1} (1977), 1501--1515.


\bibitem{V2} L. Vuillon {\em A characterization of Sturmian words by
  return words,} to appear in European. J. of Combin..

\end{thebibliography}
 \end{document}